# NUMBER OF PATHS VERSUS NUMBER OF BASIS FUNCTIONS IN AMERICAN OPTION PRICING[1]

By Paul Glasserman and Bin Yu

*Columbia University*

An American option grants the holder the right to select the time at which to exercise the option, so pricing an American option entails solving an optimal stopping problem. Difficulties in applying standard numerical methods to complex pricing problems have motivated the development of techniques that combine Monte Carlo simulation with dynamic programming. One class of methods approximates the option value at each time using a linear combination of basis functions, and combines Monte Carlo with backward induction to estimate optimal coefficients in each approximation. We analyze the convergence of such a method as both the number of basis functions and the number of simulated paths increase. We get explicit results when the basis functions are polynomials and the underlying process is either Brownian motion or geometric Brownian motion. We show that the number of paths required for worst-case convergence grows exponentially in the degree of the approximating polynomials in the case of Brownian motion and faster in the case of geometric Brownian motion.

**1. Introduction.** An American option grants the holder the right to select the time at which to exercise the option, and in this differs from a European option which may be exercised only at a fixed date. A standard result in the theory of contingent claims states that the equilibrium price of an American option is its value under an optimal exercise policy (see, e.g., Chapter 8 of [6]). Pricing an American option thus entails solving an optimal stopping problem, typically with a finite horizon.

Solving this optimal stopping problem and pricing an American option are relatively straightforward in low dimensions. Assuming a Markovian formulation of the problem, the relevant dimension is the dimension of the state

Received June 2003; revised November 2003.

[1]Supported in part by NSF Grant DMS-00-74637.

*AMS 2000 subject classifications.* Primary 60G40; secondary 65C05, 65C50, 60G35.

*Key words and phrases.* Optimal stopping, Monte Carlo methods, dynamic programming, orthogonal polynomials, finance.

vector, and this is ordinarily at least as large as the number of underlying assets on which the payoff of the option depends. In up to about three dimensions, the problem can be solved using a variety of numerical methods, including binomial lattices, finite-difference methods and techniques based on variational inequalities. (See, e.g., Chapter 5 of [10] or Chapter 9 of [19] for an introduction to these methods.) But many problems arising in practice have much higher dimensions, and these applications have motivated the development of Monte Carlo methods for pricing American options. The optimal stopping problem embedded in the valuation of an American option makes this an unconventional and challenging problem for Monte Carlo.

One class of techniques, based primarily on proposals of Carrière [4], Longstaff and Schwartz [11] and Tsitsiklis and Van Roy [17, 18], provides approximate solutions by combining simulation, regression and a dynamic programming formulation of the problem. Related methods have been used to solve dynamic programming problems in other contexts; Bertsekas and Tsitsiklis [1] discuss several techniques and applications. In this approach to American option pricing, the value function describing the option price at each time as a function of the underlying state is approximated by a linear combination of basis functions; the coefficients in this representation are estimated by applying regression to the simulated paths. Such an approximation is computed at each step in a dynamic programming procedure that starts with the option value at expiration and works backward to find the value at the current time. Any such method clearly restricts the number of possible exercise dates to be finite; these dates may be specified in the terms of the option, or they may serve as a discrete-time approximation for a continuously-exercisable option.

The convergence results available to date for these methods are based on letting the number of simulated paths increase while holding the number of basis functions fixed. Tsitsiklis and Van Roy [18] prove such a result for their method and Clément, Lamberton and Protter [5] do this for the method of Longstaff and Schwartz [11]. (The two methods differ in the backward induction procedure they use to solve the dynamic programming problem.) The convergence established by these results is therefore convergence to the approximation that would be obtained if the calculations could be carried out exactly, without the sampling error associated with Monte Carlo. Convergence to the correct option price requires a separate passage to the limit in which the number of basis functions increases.

This paper considers settings in which the number of paths and number of basis functions increase together. Our objective is to determine how quickly the number of paths must grow with the number of basis functions to ensure convergence to the correct value. The growth required turns out to be surprisingly fast in the settings we analyze. We take the underlying process to be Brownian motion or geometric Brownian motion and regress against polynomials in each case. We examine conditions for convergence to hold uniformly over coefficient vectors having a fixed norm, and in this sense our results provide a type of worst-case analysis. We show that for Brownian motion, the number of polynomials $K = K_N$ for which accurate estimation is possible from $N$ paths is $O(\log N)$; for geometric Brownian motion it is $O(\sqrt{\log N})$. Thus, the number of paths must grow exponentially with the number of polynomials in the first case, faster in the second case.

Focusing on simple models allows us to give rather precise results. Our most explicit results apply to one-dimensional problems that do not require Monte Carlo methods, but we believe they are, nevertheless, relevant to higher-dimensional problems. Many high-dimensional interest rate models have dynamics that are nearly Brownian or nearly log-Brownian; see, for example, the widely used models in Chapters 14 and 15 of [13]. Our focus on polynomials helps make our results explicit and is also consistent with, for example, examples in [11] and remarks in [5]. Our analysis relies on asymptotics of moments of the functions used in the regressions. To the extent that similar asymptotics could be derived for other basis functions and underlying distributions, our approach could be used in other settings.

We prove two types of results, providing upper and lower bounds on $K$ and thus corresponding to negative and positive results, respectively. For an upper bound on $K$, it suffices to exhibit a problem for which convergence fails. For this part of the analysis we therefore consider a single-period problem—a single regression and a single step in the backward induction. The fact that an exponentially growing sample size is necessary even in a one-dimensional, single-period problem makes the result all the more compelling. For the positive results we consider an arbitrary but fixed number of steps, corresponding to a finite set of exercise opportunities. We prove a general error bound that relies on few assumptions about the underlying Markov process or basis functions, and then specialize to the case of polynomials with Brownian motion and geometric Brownian motion.

Section 2 formulates the American option pricing problem, discusses approximate dynamic programming and presents the algorithm we analyze. Section 3 undertakes the single-period analysis, first in a normal setting then in a lognormal setting. Section 4 presents results for the multiperiod case. Proofs of some of the results in Sections 3 and 4 are deferred to Sections 5 and 6, respectively.

**2. Problem formulation.** In this section we first give a general description of the American option pricing problem, then discuss approximate dynamic programming procedures and then detail the algorithm we analyze.

2.1. *The optimal stopping problem.* A general class of American option pricing problems can be formulated through an $\Re^d$-valued Markov process $\{S(t), 0 \leq t \leq T\}$, [with $S(0)$ fixed], that records all relevant financial information, including the prices of underlying assets. We restrict attention to options admitting a finite set of exercise opportunities $0 = t_0 < t_1 < t_2 < \cdots < t_m \leq T$, sometimes called Bermudan options. (We preserve the continuous-time specification of $S$ because the lengths of the intervals $t_{i+1} - t_i$ appear in some of our results.) If exercised at time $t_n$, $n = 0, 1, \ldots, m$, the option pays $h_n(S(t_n))$, for some known functions $h_0, h_1, \ldots, h_m$ mapping $\Re^d$ into $[0, \infty)$. Let $\mathcal{T}_n$ denote the set of stopping times (with respect to the history of $S$) taking values in $\{t_n, t_{n+1}, \ldots, t_m\}$ and define

$$(1) \qquad V_n^*(x) = \sup_{\tau \in \mathcal{T}_n} \mathsf{E}[h_\tau(S(\tau))|S(t_n) = x], \qquad x \in \Re^d,$$

for $n = 0, 1, \ldots, m$. Then $V_n^*(x)$ is the value of the option at $t_n$ in state $x$, given that the option was not exercised prior to $t_n$. For simplicity, we have not included explicit discounting in (1). Deterministic discounting can be absorbed into the definition of the functions $h_n$, and stochastic discounting can usually be accommodated in this formulation at the expense of increasing the dimension of $S$.

The option values satisfy the dynamic programming equations

$$(2) \qquad V_m^*(x) = h_m(x),$$

$$(3) \qquad V_n^*(x) = \max\{h_n(x), \mathsf{E}[V_{n+1}^*(S(t_{n+1}))|S(t_n) = x]\},$$

$n = 0, 1, \ldots, m - 1$. These can be rewritten in terms of continuation values

$$C_n^*(x) = \mathsf{E}[V_{n+1}^*(S(t_{n+1}))|S(t_n) = x], \qquad n = 0, 1, \ldots, m-1,$$

as

$$(4) \qquad C_m^*(x) = 0,$$

$$(5) \qquad C_n^*(x) = \mathsf{E}[\max\{h_{n+1}(S(t_{n+1})), C_{n+1}^*(S(t_{n+1}))\}|S(t_n) = x],$$

$n = 0, 1, \ldots, m - 1$. The option values satisfy

$$V_n^*(x) = \max\{h_n(x), C_n^*(x)\},$$

so these can be calculated from the continuation values.

2.2. *Approximate dynamic programming.* Exact calculation of (2)–(3) or (4)–(5) is often impractical, and even estimation by Monte Carlo is challenging because of the difficulty of estimating the conditional expectations in these equations. Approximate dynamic programming procedures replace these conditional expectations with linear combinations of known functions, sometimes called "features" but more commonly referred to as basis functions. Thus, for each $n = 1, \ldots, m$, let $\psi_{nk}$, $k = 0, \ldots, K$, be functions from $\Re^d$ to $\Re$ and consider approximations of the form

$$C_n^*(x) \approx \sum_{k=0}^{K} \beta_{nk} \psi_{nk}(x),$$

for some constants $\beta_{nk}$, or the corresponding approximation for $V_n^*$. Working with approximations of this type reduces the problem of finding the functions $C_n^*$ to one of finding the coefficients $\beta_{nk}$. The methods of Longstaff and Schwartz [11] and Tsitsiklis and Van Roy [17, 18] select coefficients through least-squares projection onto the span of the basis functions. Other methods applying Monte Carlo to solve (2)–(3) include Broadie and Glasserman [2, 3], Haugh and Kogan [9] and Rogers [16]; for an overview, see Glasserman [8].

To simplify notation, we write $S_n$ for $S(t_n)$. We write $\psi_n$ for the vector of functions $(\psi_{n0}, \ldots, \psi_{nK})^\top$. The following basic assumptions will be in force throughout:

(A0) $\psi_{n0} \equiv 1$ for $n = 1, \ldots, m$; $\mathsf{E}[\psi_n(S_n)] = 0$, for $n = 1, \ldots, m$; and

$$\Psi_n = \mathsf{E}[\psi_n(S_n) \psi_n(S_n)^\top]$$

is finite and nonsingular, $n = 1, \ldots, m$.

For any square-integrable random variable $Y$ define the projection

$$\Pi_n Y = \psi_n^\top(S_n) \Psi_n^{-1} \mathsf{E}[Y \psi_n(S_n)].$$

Thus,

(6) $$\Pi_n Y = \sum_{k=0}^{K} a_k \psi_{nk}(S_n)$$

with

(7) $$(a_0, \ldots, a_K)^\top = \Psi_n^{-1} \mathsf{E}[Y \psi_n(S_n)].$$

We also write

$$(\Pi_n Y)(x) = \sum_{k=0}^{K} a_k \psi_{nk}(x)$$

for the function defined by the coefficients (7).

Define an approximation to (4)–(5) as follows: $C_m(x) \equiv 0$,

(8) $$C_n(x) = (\Pi_n \max\{h_{n+1}(S_{n+1}), C_{n+1}(S_{n+1})\})(x).$$

As in (6), the application of the projection $\Pi_n$ results in a linear combination of the basis functions, so

$$\text{(9)} \quad C_n(x) = (\Pi_n \max\{h_{n+1}, C_{n+1}\})(x) = \sum_{k=0}^{K} \beta_{nk} \psi_{nk}(x)$$

with $\beta_n^\top = (\beta_{n0}, \ldots, \beta_{nK})$ defined as in (7) but with $Y$ replaced by

$$\text{(10)} \quad V_{n+1}(S_{n+1}) \equiv \max\{h_{n+1}(S_{n+1}), C_{n+1}(S_{n+1})\}.$$

With the payoff functions $h_n$ fixed, we can rewrite (9) using the operator

$$\text{(11)} \quad L_n C_{n+1} = \Pi_n(\max\{h_{n+1}, C_{n+1}\}).$$

Exact calculation of the projection in (8) is usually infeasible, but it is relatively easy to evaluate a sample counterpart of this recursion defined from a finite set of simulated paths of the process $S$. We consider the following procedure to approximate the coefficient vectors $\beta_n$ and the continuation values $C_n$.

*Step* 1. Set $\hat{C}_m = 0$ and $\hat{V}_m = \max\{h_m, \hat{C}_m\} = h_m$.

*Step* 2. For each $n = 1, \ldots, m-1$, repeat the following steps: Generate $N$ paths $\{S_1^{(i)}, \ldots, S_{n+1}^{(i)}\}$, $i = 1, \ldots, N$, up to time $t_{n+1}$, independent of each other and of all previously generated paths. Calculate

$$\hat{\gamma}_n = \frac{1}{N} \sum_{i=1}^{N} \hat{V}_{n+1}(S_{n+1}^{(i)}) \psi_n(S_n^{(i)}),$$

calculate the coefficients $\hat{\beta}_n = \Psi_n^{-1} \hat{\gamma}_n$ and set

$$\text{(12)} \quad \hat{C}_n = \hat{\beta}_n^\top \psi_n \equiv \hat{L}_n \hat{C}_{n+1} \equiv \hat{\Pi}_n \max\{h_{n+1}, \hat{C}_{n+1}\},$$

$$\text{(13)} \quad \hat{V}_n = \max\{h_n, \hat{C}_n\}.$$

*Step* 3. Set $\hat{C}_0(S_0) = N^{-1} \sum_{i=1}^{N} \hat{V}_1(S_1^{(i)})$ and $\hat{V}_0(S_0) = \max\{h_0(S_0), \hat{C}_0(S_0)\}$.

A few aspects of this algorithm require comment. In Step 3 we simply average the estimated values at $t_1$ to get the continuation value at time 0 because $S(0)$ is fixed. The operators $\hat{L}_n$ and $\hat{\Pi}_n$ implicitly defined in (12) are the sample counterparts of those in (6) and (11), using estimated rather than exact coefficients. The coefficient estimates in Step 2 use the matrices $\Psi_n$. In ordinary least-squares regression, each $\Psi_n$ would be replaced with its sample counterpart,

$$\frac{1}{N} \sum_{i=1}^{N} \psi_n(S_n^{(i)}) \psi_n(S_n^{(i)})^\top,$$

calculated from the simulated values themselves. (Owen [14] calls the use of the exact matrix *quasi*-regression.) In our examples, the $\Psi_n$ are indeed available explicitly and using this formulation simplifies the analysis.

In Step 2 we have used an independent set of paths to estimate coefficients at each date, though the algorithms of Longstaff and Schwartz [11] and Tsitsiklis and Van Roy [17, 18] use a single set of paths for all dates. This modification is theoretically convenient because it makes the coefficients of $\hat{C}_{n+1}$ independent of the points at which $\hat{C}_{n+1}$ is evaluated in the calculation of $\hat{\gamma}_n$. This distinction is relevant only to the multiperiod analysis of Section 4 and disappears in the single-period analysis of Section 3. The worst case over all multiperiod problems is at least as bad as the worst single-period problem. The results in Section 3 thus provide lower bounds on the worst-case convergence rate for multiperiod problems whether one uses independent paths at each date or a single set of paths for all dates.

**3. Single-period problem.** For the single-period problem, we fix dates $t_1 < t_2$ and consider the estimation of coefficients $\beta_0, \ldots, \beta_K$ in the projection of a function of $S_2$ onto the span of $\psi_{1k}(S_1)$, $k = 0, \ldots, K$. Thus,

$$\beta = (\beta_0, \ldots, \beta_K)^\top = \Psi^{-1} \gamma \tag{14}$$

with $\Psi = \Psi_1$ and $\gamma = \mathsf{E}[Y \psi_1(S_1)]$ for some $Y$. In a simplified instance of the algorithm of the previous section, we simulate $N$ independent copies $(S_1^{(i)}, Y^{(i)})$, $i = 1, \ldots, N$, and compute the estimate

$$\tilde{\beta} = \Psi^{-1} \tilde{\gamma}, \tag{15}$$

where $\tilde{\gamma}$ is the unbiased estimator of $\gamma$ with components

$$\tilde{\gamma}_k = \frac{1}{N} \sum_{i=1}^{N} Y^{(i)} \psi_{1k}(S_1^{(i)}), \qquad k = 0, 1, \ldots, K. \tag{16}$$

We analyze the convergence of $\tilde{\beta}$ (and $\tilde{\gamma}$) as both $N$ and $K$ increase.

We denote by $|x|$ the Euclidean norm of the vector $x$. For a matrix $A$, we denote by $\|A\|$ the Euclidean matrix norm, meaning the square root of the sum of squared elements of $A$. It follows that $|Ax| \leq \|A\| |x|$ and then from (14) and (15),

$$\frac{1}{\|\Psi\|} |\tilde{\gamma} - \gamma| \leq |\tilde{\beta} - \beta| \leq \|\Psi^{-1}\| |\tilde{\gamma} - \gamma|. \tag{17}$$

The Euclidean norm on vectors is a measure of the proximity of the functions determined by vectors of coefficients. To make this more explicit, let $b$ and $c$ be coefficient vectors and let $S_n$ have density $g_n$. Then

$$\int \left( \sum_{k=0}^{K} b_k \psi_{nk}(x) - \sum_{k=0}^{K} c_k \psi_{nk}(x) \right)^2 g_n(x) \, dx = (b-c)^\top \Psi_n (b-c)$$

and

$$\frac{1}{\|\Psi_n^{-1}\|}|b-c|^2 \leq (b-c)^\top \Psi_n (b-c) \leq \|\Psi_n\|\,|b-c|^2.$$

Thus, the Euclidean norm on vectors gives the $L^2$ norm (with respect to $g_n$) for the functions determined by the vectors, up to factors of $\|\Psi_n\|$ and $\|\Psi_n^{-1}\|$ that will prove to be negligible in the settings we consider.

We therefore investigate the convergence of the expected squared difference $\mathsf{E}[|\beta - \tilde{\beta}|^2]$. Because this is the mean square error of $\tilde{\beta}$, we also denote it by $\mathsf{MSE}(\tilde{\beta})$. Thus, (17) implies

(18) $$\frac{1}{\|\Psi\|^2}\mathsf{E}[|\tilde{\gamma}-\gamma|^2] \leq \mathsf{MSE}(\tilde{\beta}) \leq \|\Psi^{-1}\|^2 \mathsf{E}[|\tilde{\gamma}-\gamma|^2].$$

For a given number of replications $N$ and basis functions $K$, $\mathsf{MSE}(\tilde{\beta})$ can be made arbitrarily large or small by multiplying $\beta$ by a constant. To get meaningful results, we therefore adopt the following normalization:

(A1) $|\beta| = 1$.

We investigate the convergence of the supremum of the $\mathsf{MSE}(\tilde{\beta})$ over all $\beta$ satisfying this condition. In order to investigate how $N$ must grow with $K$, we assume that the regression representation is, in fact, valid, in a sense implied by the following two conditions:

(A2) $Y$ has the form

$$Y = \sum_{k=0}^{K} a_k \psi_{2k}(S_2),$$

for some constants $a_k$.

(A3) There exist functions $f_k : \Re_+ \to \Re_+$, $k = 0, \ldots, K$, such that

$$\mathsf{E}[f_k(t_2)\psi_{2k}(S_2)|S_1] = f_k(t_1)\psi_{1k}(S_1), \qquad t_2 \geq t_1.$$

Condition (A3) states that the $\psi_{nk}(S_n)$ are martingales, up to a deterministic function of time. Condition (A2), though a strong assumption in practice, makes Theorems 1 and 2 more compelling: the rapid growth in the number of paths implied by the theorems holds even though we have chosen the "correct" basis functions, in the sense of (A2). The results of Section 4 give sufficient conditions for convergence without such an assumption.

Under assumptions (A2) and (A3), we have

$$\gamma_k = \mathsf{E}[Y\psi_{1k}(S_1)]$$

(19) $$= \mathsf{E}\left[\sum_{l=0}^{K} a_l \psi_{2l}(S_2)\psi_{1k}(S_1)\right]$$

$$= \sum_{l=0}^{K} a_l \frac{f_l(t_1)}{f_l(t_2)} \mathsf{E}[\psi_{1l}(S_1)\psi_{1k}(S_1)].$$

The restriction on $\beta$ in (A1) then restricts $a$.

Returning to the analysis of $\mathsf{MSE}(\tilde{\beta})$, (18) indicates that we need to analyze the mean square error of $\tilde{\gamma}$, for which (since $\mathsf{E}[\tilde{\gamma}] = \gamma$) we get

$$\mathsf{E}[|\tilde{\gamma} - \gamma|^2] = \sum_{k=0}^{K} \mathsf{Var}[\tilde{\gamma}_k] \tag{20}$$

$$= \sum_{k=0}^{K} \frac{1}{N} \mathsf{E}[Y^2 \psi_{1k}^2(S_1)] - \frac{1}{N} \sum_{k=0}^{K} \gamma_k^2. \tag{21}$$

Thus, using (18), (A2) and the Cauchy–Schwarz inequality,

$$\mathsf{MSE}(\tilde{\beta}) \leq \|\Psi^{-1}\|^2 \mathsf{E}[|\tilde{\gamma} - \gamma|^2]$$

$$\leq \|\Psi^{-1}\|^2 \sum_{k=0}^{K} \frac{1}{N} \mathsf{E}[Y^2 \psi_{1k}^2(S_1)] \tag{22}$$

$$\leq \|\Psi^{-1}\|^2 \frac{1}{N} \sum_{l=0}^{K} a_l^2 \sum_{k,j=0}^{K} \mathsf{E}[\psi_{2j}^2(S_2) \psi_{1k}^2(S_1)].$$

To get a lower bound, we may define $Y^* = a_K^* \psi_{2K}(S_2)$, with $a_K^*$ chosen such that the corresponding $\beta^*$ satisfies $|\beta^*| = 1$. Using (18) and (20), we then get

$$\sup_{|\beta|=1} \mathsf{MSE}(\tilde{\beta}) \geq \frac{1}{\|\Psi\|^2} \left( \sum_{k=0}^{K} \frac{1}{N} \mathsf{E}[Y^{*2} \psi_{1k}^2(S_1)] - \frac{1}{N} \sum_{k=0}^{K} \gamma_k^2 \right)$$

$$= \frac{1}{\|\Psi\|^2} \left( a_K^{*2} \sum_{k=0}^{K} \frac{1}{N} \mathsf{E}[\psi_{2K}^2(S_2) \psi_{1k}^2(S_1)] - \frac{1}{N} \sum_{k=0}^{K} \gamma_k^2 \right). \tag{23}$$

From (22) and (23) we see that the key to the analysis of the uniform convergence of $\mathsf{MSE}(\tilde{\beta})$ lies in the growth of fourth-order moments of the form $\mathsf{E}[\psi_{2j}^2(S_2)\psi_{1k}^2(S_1)]$. This, in turn, depends on the choice of basis functions and on the law of the underlying process $S$. We analyze the case of polynomials with Brownian motion and geometric Brownian motion.

3.1. *Normal setting.* For this section, let $\{S(t), 0 \leq t \leq T\}$ be a standard Brownian motion. We define the basis functions through the Hermite polynomials

$$H_{e_n}(x) = \sum_{i=0}^{\lfloor n/2 \rfloor} \frac{(-1)^i n! x^{n-2i}}{(n-2i)! i! 2^i}, \qquad n = 0, 1, \ldots,$$

where $\lfloor n/2 \rfloor$ denotes the integer part of $n/2$. The Hermite polynomials have the following useful properties: They are orthogonal with respect to the standard normal density $\phi$, in the sense that $H_{e_0} \equiv 1$ and

$$\int H_{e_i}(x) H_{e_j}(x) \phi(x)\, dx = \begin{cases} 0, & i \neq j, \\ i!, & i = j. \end{cases}$$

They define martingales, in the sense that (see, e.g., [15], page 151)

$$\mathsf{E}\left[ t_2^{i/2} H_{e_i}\left(\frac{S(t_2)}{\sqrt{t_2}}\right) \Big| S(t_1) \right] = t_1^{i/2} H_{e_i}\left(\frac{S(t_1)}{\sqrt{t_1}}\right),$$

for $t_2 \geq t_1$. And their squares admit the expansion

(24) $$(H_{e_n}(x))^2 = (n!)^2 \sum_{i=0}^{n} \frac{H_{e_{2i}}(x)}{(i!)^2 (n-i)!}.$$

The functions

(25) $$\psi_{nk}(x) = \frac{1}{\sqrt{k!}} H_{e_k}(x/\sqrt{t_n})$$

satisfy (A3) with $f_k(t) = t^{k/2}$. They are also orthogonal and their $\Psi$ matrix is the identity. Thus, $\beta = \gamma$ and $\tilde{\beta} = \tilde{\gamma}$.

We can now state the main result of this section. Let $\rho = t_2/t_1$, and for $\rho \geq 1$ define

$$c_\rho = 2 \log(2 + \sqrt{\rho}).$$

TABLE 1
*Estimates of* $\mathsf{MSE}(\tilde{\beta})$ *for various combinations of $K$ basis functions and $N$ paths. The critical values $K = \log N / c_\rho$ are displayed by in the bottom row and also indicated by the horizontal line through the table*

| | | | | | $N$ | | | | |
|---|---|---|---|---|---|---|---|---|---|
| $K$ | 500 | 1000 | 2000 | 4000 | 8000 | 16000 | 32000 | 64000 | 128000 |
| 1 | 0.01 | 0.00 | 0.00 | 0.00 | 0.00 | 0.00 | 0.00 | 0.00 | 0.00 |
| 2 | 0.08 | 0.04 | 0.02 | 0.01 | 0.00 | 0.00 | 0.00 | 0.00 | 0.00 |
| 3 | 0.67 | 0.31 | 0.17 | 0.08 | 0.04 | 0.02 | 0.01 | 0.00 | 0.00 |
| 4 | 5.6 | 3.0 | 1.6 | 0.73 | 0.36 | 0.18 | 0.09 | 0.05 | 0.02 |
| 5 | 52.7 | 23.4 | 13.5 | 6.0 | 3.1 | 1.5 | 0.8 | 0.40 | 0.20 |
| 6 | 427.2 | 155.7 | 93.3 | 38.4 | 24.0 | 10.8 | 6.2 | 3.1 | 1.5 |
| 7 | 2403 | 1202 | 600.8 | 300.4 | 150.2 | 75.1 | 37.5 | 18.8 | 9.4 |
| 8 | 11447 | 5723 | 2862 | 1431 | 715.4 | 357.7 | 178.9 | 89.4 | 44.7 |
| 9 | | | 9856 | 4928 | 2464 | 1232 | 616 | 308 | 154 |
| 10 | | | | | 6109 | 3054 | 1527 | 764 | 381 |
| 11 | | | | | | | 2810 | 1405 | 702 |
| 12 | | | | | | | | | 1023 |
| Bound | 2.5 | 2.8 | 3.1 | 3.4 | 3.7 | 3.9 | 4.2 | 4.5 | 4.8 |

THEOREM 1. *Let $\psi_{nk}$ be as in (25) and suppose* (A2) *holds. If $K = (1-\delta) \times \log N / c_\rho$ for some $\delta > 0$, then*

$$\lim_{N \to \infty} \sup_{|\beta|=1} \mathsf{MSE}(\tilde{\beta}) = 0. \tag{26}$$

*If $K = (1+\delta) \log N / c_\rho$ for some $\delta > 0$, then*

$$\lim_{N \to \infty} \sup_{|\beta|=1} \mathsf{MSE}(\tilde{\beta}) = \infty. \tag{27}$$

This result shows rather precisely that, from a sample size of $N$, the highest $K$ for which coefficients of polynomials of order $K$ can be estimated uniformly well is $O(\log N)$. Equivalently, the sample size required to achieve convergence grows exponentially in $K$.

This is illustrated numerically in Table 1, which shows estimates of $\mathsf{MSE}(\tilde{\beta})$ for various combinations of $N$ and $K$. The results shown are for $Y = \rho^{K/2} H_{e_K}(S_2/\sqrt{t_2})/\sqrt{K!}$, with $t_1 = 1$ and $t_2 = 2$, a special case of the $Y$ we use to prove (27). The estimates are computed as follows. For each entry of the table, we generate 5000 batches, each consisting of $N$ paths. From each batch we compute $\tilde{\beta}$ and then take the average of $|\tilde{\beta} - \beta|^2$ over the 5000 batches. This average provides our estimate of $\mathsf{MSE}(\tilde{\beta})$ in each case with $K \leq 6$. For $K \geq 7$ this produced unacceptably high variability, so for those cases we calculated $\mathsf{MSE}(\tilde{\beta})$ from 5000 replications of $N = 500{,}000$ and then scaled the estimate by $N$.

The bottom row of the table displays the critical values $K = \log N / c_\rho$ provided by Theorem 1; these values are also indicated by the horizontal line through the table. As indicated by the theorem, $\mathsf{MSE}(\tilde{\beta})$ explodes along any diagonal line through the table steeper than the critical line, and remains small above the critical line.

The proof of the theorem uses the following two lemmas, proved in Section 5.

LEMMA 1. *For the $\psi_{nk}$ in (25) and $\rho = t_2/t_1$,*

$$\mathsf{E}[\psi_{2k_2}(S_2)\psi_{1k_1}(S_1)] = \begin{cases} 0, & k_1 \neq k_2, \\ \rho^{-k_1/2}, & k_1 = k_2, \end{cases} \tag{28}$$

$$\mathsf{E}[\psi_{2k_2}(S_2)^2 \psi_{1k_1}(S_1)^2] = \sum_{k=0}^{k_1 \wedge k_2} \rho^{-k} \binom{2k}{k} \binom{k_1}{k} \binom{k_2}{k}, \tag{29}$$

*with $k_1 \wedge k_2$ the minimum of $k_1$ and $k_2$. Equation* (29) *is strictly increasing in $k_1$ and $k_2$.*

For the special case $k_1 = k_2 = K$, (29) yields

$$\text{E}[\psi_{2K}(S_2)^2\psi_{1K}(S_1)^2] = \sum_{k=0}^{K}\rho^{-k}\binom{2k}{k}\binom{2K}{k}^2. \tag{30}$$

As a step toward bounding this expression, let $k^*$ denote the index of the largest summand so that

$$\rho^{-k^*}\binom{2k^*}{k^*}\binom{K}{k^*}^2 = \max_{0\leq k\leq K}\rho^{-k}\binom{2k}{k}\binom{K}{k}^2. \tag{31}$$

For $k^*$, we have the following lemma.

LEMMA 2. *As $K \to \infty$,*

$$k^* = \frac{2}{2+\sqrt{\rho}}K(1+o(1)).$$

PROOF OF THEOREM 1. We bound $\text{MSE}(\tilde{\beta})$ from above based on (22). Combining the fact that $\beta = \gamma$ (because $\Psi = I$) with (19) and (28) we get

$$\beta_k = \sum_{l=0}^{K} a_l \text{E}[\psi_{2l}(S_2)\psi_{1k}(S_1)] = a_k\rho^{-k/2}.$$

Thus, $|\beta|=1$ implies $a_k^2 \leq \rho^k$. From (30) and (31), we get

$$\rho^{-k^*}\binom{2k^*}{k^*}\binom{K}{k^*}^2 < \text{E}[\psi_{2K}(S_2)^2\psi_{1K}(S_1)^2] \tag{32}$$

$$< (K+1)\rho^{-k^*}\binom{2k^*}{k^*}\binom{K}{k^*}^2.$$

Recalling that $\Psi = I$ and applying the inequality $a_k^2 \leq \rho^k$ to (22) we get

$$\sup_{|\beta|=1}\text{MSE}(\tilde{\beta}) \leq \sup_{|\beta|=1}\|\Psi^{-1}\|^2\frac{1}{N}\sum_{l=0}^{K}a_l^2\sum_{k,l=0}^{K}\text{E}[\psi_{2l}^2(S_2)\psi_{1k}^2(S_1)]$$

$$\leq (K+1)\frac{1}{N}\sum_{k=0}^{K}\rho^k\sum_{k,l=0}^{K}\text{E}[\psi_{2k}^2(S_2)\psi_{1l}^2(S_1)]$$

$$< \frac{(K+1)^2}{N}\rho^K(K+1)^2\text{E}[\psi_{2K}^2(S_2)\psi_{1K}^2(S_1)] \tag{33}$$

$$< \frac{(K+1)^5}{N}\rho^{K-k^*}\binom{2k^*}{k^*}\binom{K}{k^*}^2, \tag{34}$$

where (33) follows from Lemma 1 and the last inequality follows from (32).

To get a lower bound on the supremum of $\mathsf{MSE}(\tilde{\beta})$ we use (23) with

$$Y^* = \rho^{K/2}\psi_{2K}(S_2) \equiv a_K^* \psi_{2K}(S_2),$$

for which $\beta_K = 1$ and $\beta_k = 0$, $k \neq K$. By applying Lemma 1 and the lower bound in (32), (23) becomes

$$\sup_{|\beta|=1} \mathsf{MSE}(\tilde{\beta}) \geq \frac{1}{\|\Psi\|^2} \frac{1}{N} \left( a_K^{*2} \sum_{k=0}^{K} \mathsf{E}[\psi_{2K}^2(S_2)\psi_{1k}^2(S_1)] - 1 \right)$$

$$\geq \frac{1}{\|\Psi\|^2} \frac{a_K^{*2}}{N} \mathsf{E}[\psi_{2K}^2(S_2)\psi_{1K}^2(S_1)]$$

(35) $$\geq \frac{1}{K+1}\frac{1}{N}\rho^{K-k^*}\binom{2k^*}{k^*}\binom{K}{k^*}^2.$$

By Stirling's approximation $n! \sim \sqrt{2n\pi}(\frac{n}{e})^n$ and Lemma 2 we get

$$\binom{K}{k^*} = \frac{K!}{k^*!(K-k^*)!}$$

$$= \frac{\sqrt{2K\pi}(K/e)^K(1+o(1))}{\sqrt{2k^*\pi}\sqrt{2(K-k^*)\pi}(k^*/e)^{k^*}(K-k^*/e)^{K-k^*}(1+o(1))}$$

(36) $$= \frac{1}{\sqrt{2abK\pi}a^{aK}b^{bK}}(1+o(1)),$$

with $a = 2/(2+\sqrt{\rho})$ and $b = 1-a$. Also,

$$\binom{2k^*}{k^*} = \frac{2k^*!}{k^*!k^*!}$$

$$= \frac{\sqrt{4k^*\pi}(2k^*/e)^{2k^*}(1+o(1))}{2k^*\pi(k^*/e)^{2k^*}(1+o(1))}$$

(37) $$= \frac{2^{2aK}}{\sqrt{aK\pi}}(1+o(1)).$$

By substituting (36) and (37) into (34) and (35) we get

$$\frac{\rho^{bK}2^{2aK}}{2N(K+1)\sqrt{aK\pi}abK\pi a^{2aK}b^{2bK}}(1+o(1))$$

$$\leq \sup_{|\beta|=1} \mathsf{MSE}(\tilde{\beta})$$

$$\leq \frac{(K+1)^5\rho^{bK}2^{2aK}}{2N\sqrt{aK\pi}abK\pi a^{2aK}b^{2bK}}(1+o(1)).$$

Simple algebra verifies that $c_\rho = 2a\log(2) - 2a\log(a) - 2b\log(b) + b\log(\rho)$, so we can rewrite these bounds as

$$\frac{e^{c_\rho K}}{2N(K+1)\sqrt{aK\pi abK\pi}}(1+o(1)) \leq \sup_{|\beta|=1} \mathsf{MSE}(\tilde{\beta})$$

$$\leq \frac{(K+1)^5 e^{c_\rho K}}{2N\sqrt{aK\pi abK\pi}}(1+o(1)).$$

If $K = (1-\delta)\log N/c_\rho$ for some $\delta > 0$, then as $N \to \infty$,

$$\log\left\{\frac{(K+1)^5 e^{c_\rho K}}{2N\sqrt{aK\pi abK\pi}}(1+o(1))\right\} = -\delta \log N + o(\log N) \to -\infty,$$

so (26) holds. If $K = (1+\delta)\log N/c_\rho$ for some $\delta > 0$, then as $N \to \infty$,

$$\log\left\{\frac{e^{c_\rho K}}{2N(K+1)\sqrt{aK\pi abK\pi}}(1+o(1))\right\} = \delta \log N + o(\log N) \to \infty,$$

and (27) holds. □

3.2. *Lognormal setting.* We now take $S$ to be geometric Brownian motion, $S(t) = \exp(W(t) - t/2)$, with $W$ a standard Brownian motion. For the basis functions $\psi_{nk} \equiv \psi_k$, we use multiples of the powers $x^k$ to get the martingales

(38) $$\psi_k(S(t)) = e^{kW(t) - k^2 t/2}.$$

These functions satisfy (A0). The main result of this section is the following:

THEOREM 2. *Let the $\psi_k$ be as in (38) and suppose (A2) holds. If*

$$K = \sqrt{\frac{(1-\delta)\log N}{5t_1 + t_2}}$$

*for some $\delta > 0$, then*

$$\lim_{N\to\infty} \sup_{|\beta|=1} \mathsf{MSE}(\tilde{\beta}) = 0.$$

*If*

$$K = \sqrt{\frac{(1+\delta)\log N}{3t_1 + t_2}}$$

*for some $\delta > 0$, then*

$$\lim_{N\to\infty} \sup_{|\beta|=1} \mathsf{MSE}(\tilde{\beta}) = \infty.$$

Compared with the normal case in Theorem 1, we see that here $K$ must be much smaller—of the order of $\sqrt{\log N}$. Accordingly, $N$ must be much larger—of the order of $\exp(K^2)$. The analysis in this setting is somewhat more complicated than in the normal case because the $\psi_k$ are no longer orthogonal. To prove the theorem we state some lemmas that are proved in Section 5.

LEMMA 3. *For $t_2 \geq t_1$ and $k_1, k_2 = 0, \ldots, K$,*
$$\mathsf{E}[\psi_{k_1}(S_1)\psi_{k_2}(S_2)] = e^{k_1 k_2 t_1},$$
$$\mathsf{E}[\psi_{k_1}^2(S_1)\psi_{k_2}^2(S_2)] = e^{k_1^2 t_1 + k_2^2 t_2 + 4 k_1 k_2 t_1},$$
*and $\mathsf{E}[\psi_{k_1}(S_1)^2 \psi_{k_2}(S_2)^2]$ is strictly increasing in $k_1$ and $k_2$.*

Using the first statement in the lemma, we find that the matrix $\Psi(t)$ with $ij$th entry $\mathsf{E}[\psi_{i-1}(S(t))\psi_{j-1}(S(t))]$ is given by

$$\Psi(t) = \begin{pmatrix} 1 & 1 & 1 & \cdots & 1 \\ 1 & e^t & e^{2t} & \cdots & e^{Kt} \\ 1 & e^{2t} & e^{4t} & \cdots & e^{2Kt} \\ \vdots & \vdots & \vdots & \ddots & \vdots \\ 1 & e^{Kt} & e^{2Kt} & \cdots & e^{K^2 t} \end{pmatrix}.$$

We write $\Psi$ for $\Psi(t_1)$.

LEMMA 4. *We have $\|\Psi(t)\| \leq (K+1)^2 e^{2K^2 t}$ and, with $C(t) = \exp(-2e/(e^t - 1)^2)$,*
$$\|\Psi(t)^{-1}\| \leq C^{-1}(t) K(K+1) \left(\frac{e^t}{e^t - 1}\right)^K.$$

PROOF OF THEOREM 2. Condition (A2) and the martingale property of the $\psi_k(S(t))$ imply that
$$\mathsf{E}[Y|S_1] = \sum_{k=0}^{K} a_k \psi_k(S_1),$$

and, thus, that $\beta_k = a_k$, $k = 0, 1, \ldots, K$. In this case, the normalization $|\beta| = 1$ is equivalent to $|(a_0, \ldots, a_K)| = 1$. Applying this in (22) and then applying Lemmas 3 and 4 we get

$$\sup_{|\beta|=1} \mathsf{MSE}(\tilde{\beta}) \leq \sup_{|\beta|=1} \|\Psi^{-1}\|^2 \frac{1}{N} \mathsf{E}\left[\sum_{k=0}^{K} \psi_k^2(S_2) \psi_K^2(S_1)\right]$$

$$\leq \|\Psi^{-1}\|^2 \frac{1}{N}(K+1)\mathsf{E}[\psi_K^2(S_2)\psi_K^2(S_1)]$$

$$\leq C^{-2}(t_1)K^2(K+1)^2\left(\frac{e^{t_1}}{e^{t_1}-1}\right)^{2K}\frac{K+1}{N}e^{5K^2t_1+K^2t_2}.$$

If we now take $K = \sqrt{\frac{(1-\delta)\log N}{5t_1+t_2}}$, then as $N \to \infty$,

$$\log\left\{C(t_1)^2 K^2(K+1)^3\left(\frac{e^{t_1}}{e^{t_1}-1}\right)^{2K}\frac{1}{N}e^{5K^2t_1+K^2t_2}(1+o(1))\right\}$$

$$= -\delta \log N + o(\log N) \to -\infty,$$

which proves the first assertion in the theorem.

For the second part of the theorem, define

$$Y^* = e^{KW(t_2)-K^2t_2/2},$$

for which $\beta^*$ is $(0,\ldots,0,1)^\top$. The corresponding vector $\gamma^*$ is $\Psi\beta^*$, the last column of $\Psi$. Applying this in (23) and using Lemmas 3 and 4 we get

$$\sup_{|\beta|=1} \mathsf{MSE}(\tilde{\beta}) \geq \frac{1}{\|\Psi\|^2}\left(\sum_{k=0}^{K}\frac{1}{N}\mathsf{E}[\psi_K^2(S_2)\psi_k^2(S_1)] - \frac{1}{N}\sum_{k=0}^{K}(\gamma_k^*)^2\right)$$

$$\geq \frac{1}{\|\Psi\|^2}\frac{1}{N}(\mathsf{E}[\psi_K^2(S_2)\psi_K^2(S_1)] - (\gamma_K^*)^2)$$

$$\geq \frac{1}{N(K+1)^2 e^{2K^2t_1}}(e^{5K^2t_1+K^2t_2} - e^{2K^2t_1})$$

$$= \frac{1}{N(K+1)^2}e^{3K^2t_1+K^2t_2}(1+o(1)).$$

If we now take $K = \sqrt{\frac{(1+\delta)\log N}{3t_1+t_2}}$, then as $N \to \infty$,

$$\log\left\{\frac{1}{N(K+1)^2}e^{3K^2t_1+K^2t_2}(1+o(1))\right\} = \delta \log N + o(\log N) \to \infty,$$

proving the second assertion in the theorem. $\square$

The analysis of this section differs from the normal setting of Section 3.1 in that the polynomials (38) are not orthogonal. In the Brownian case, the Hermite polynomials are orthogonal and (after appropriate scaling) martingales. In using (38), we have chosen to preserve the martingale property rather than orthogonality. As a consequence $\|\Psi^{-1}\|$ and $1/\|\Psi\|$ appear in our bounds on $\mathsf{MSE}(\tilde{\beta})$. From Lemma 4 we see that $\|\Psi^{-1}\|$ has an asymptotically negligible effect on the upper bound for $\mathsf{MSE}(\tilde{\beta})$, and with or without the factor of $1/\|\Psi\|$, the lower bound on $\mathsf{MSE}(\tilde{\beta})$ is exponential in a multiple of $K^2$. The slower convergence rate in the lognormal setting therefore does not appear to result from the lack of orthogonality.

**4. Multiperiod problem.** We now turn to conditions that ensure convergence of the multiperiod algorithm in Section 2.2 as both the number of basis functions $K$ and the number of paths $N$ increase. We first formulate a general result bounding the error in the estimated continuation values, then specialize to the normal and lognormal settings.

4.1. *General bound.* We use the following conditions.

(B1) $\mathsf{E}[\psi_{nk}^2(S_n)]$ and $\mathsf{E}[\psi_{nk}^4(S_n)]$ are increasing in $n$ and $k$.

As explained in the discussion of the single-period problem, we need some normalization on the regression coefficients in order to make meaningful statements about worst-case convergence. For a problem with $m$ exercise opportunities, we impose

(B2) $|\beta_{m-1}| = 1$.

This condition is analogous to the one we used in the single-period problem, where $\beta$ was a vector of coefficients at time $t_1$ and $Y$ was a linear combination of functions evaluated at $S(t_2)$.

We also need a condition on the functions $h_n$ that determine the payoff upon exercise at time $t_n$. The following condition turns out to be convenient:

(B3) $\mathsf{E}[h_n^4(S_n)] \leq (\frac{t_n}{t_{n-1}})^{2K} \mathsf{E}[\psi_{nK}^4(S_n)]$, for $n = 0, 1, \ldots, m$.

Suppose $S_n$ has density $g_n$ and define the weighted $L^2$ norm on functions $G : \Re \to \Re$,

$$\|G\|_n = \sqrt{\int G(x)^2 \, g_n(x) \, dx}.$$

With $\hat{C}_n$ the estimated continuation value defined by (12), we analyze the error $\mathsf{E}[\|\hat{C}_n - C_n\|_n^2]$.

We need some additional notation. Let

$$c = \max_{n=1,\ldots,m-1} \frac{t_{n+1}}{t_n}, \qquad B_K = \max_{n=1,\ldots,m-1} \|\Psi_n^{-1}\|,$$

$$H_K = \max\{c^K, B_K^2(K+1)\}, \qquad A_K = (K+1)H_K \mathsf{E}[\psi_{mK}^4(S_m)].$$

Under (A0), $B_K$ is well defined. We can now state the main result of this section.

THEOREM 3. *If assumptions* (A0) *and* (B1)–(B3) *hold, then*

$$(39) \quad \mathsf{E}[\|\hat{C}_n - C_n\|_n^2] \leq (2^{m-n} - 1)\frac{(K+1)^2}{N} B_K A_K^{m-n} (\mathsf{E}[\psi_{mK}^2(S_m)])^2 (1 + o(1)).$$

This result is proved in Section 6. Its consequences will be clearer once we illustrate it in the normal and lognormal settings.

4.2. *Multiperiod examples.*

4.2.1. *Normal setting.* As in Section 3.1, let $S$ be a standard Brownian motion and let the $\psi_{nk}$ be as in (25). Each $\Psi_n$ is then the identity matrix, $n = 1, \ldots, m$. It follows that

$$(40) \qquad B_K = \max_n \|\Psi_n^{-1}\| = \sqrt{K+1}.$$

Also,

$$H_K = \max\{c^K, B_K^2(K+1)\} = c^K$$

for all sufficiently large $K$.

To bound $\mathsf{E}[\psi_{mK}^4(S_m)]$ (which appears in $A_K$), we use (29) (with $t_1 = t_2$ and $k_1 = k_2 = K$) and then Stirling's formula and Lemma 2 to get

$$\mathsf{E}[\psi_{mK}^4(S_m)] = \sum_{k=0}^{K} \binom{2k}{k} \binom{K}{k}^2 \leq (K+1) \frac{9\sqrt{3}}{4\sqrt{2K^3\pi^3}} 3^{2K}(1+o(1)).$$

The expression on the right follows from (32), (36) and (37) upon noting that with $\rho = 1$ we get $a = 2/3$ and $b = 1/3$. Substituting this expression and (40) into (39) yields

$$\mathsf{E}[\|\hat{C}_n - C_n\|^2]$$
$$< (2^{m-n} - 1)\frac{(K+1)^{2m-2n+5/2}}{N} \left(\frac{9\sqrt{3}}{4\sqrt{2K^3\pi^3}} 3^{2K}\right)^{m-n} c^{(m-n)K}(1+o(1)).$$

It now follows that if

$$K = \frac{(1-\delta)\log N}{(m-n)(2\log 3 + \log c)}$$

for some $\delta > 0$, then

$$(41) \qquad \lim_{N \to \infty} \sup_{|\beta_{m-1}|=1} \mathsf{E}[\|\hat{C}_n - C_n\|_n^2] = 0.$$

In other words, we have convergence of the estimated continuation values at all exercise opportunities, as both $N$ and $K$ increase. If the basis functions eventually span the true optimum, in the sense that $\|C_n - C_n^*\|_n \to 0$ as $K \to \infty$, then by the triangle inequality, (41) holds with $C_n$ replaced by $C_n^*$.

On the other hand, from Theorem 1 we know that if $K = (1+\delta)\log N/c_\rho$ for any $\delta > 0$, with $\rho = t_m/t_{m-1}$, then

$$\lim_{N \to \infty} \sup_{|\beta_{m-1}|=1} \mathsf{E}[\|\hat{C}_{m-1} - C_{m-1}\|_{m-1}^2] = \infty.$$

Thus, the crititcal rate of $K$ for the multiperiod problem is $O(\log N)$, just as in the single-period problem.

4.2.2. *Lognormal setting.* Now we take $S$ to be geometric Brownian motion and use the basis functions of Section 3.2. In this case we have

$$B_K = \max_n \|\Psi_n^{-1}\| < \max_n C^{-1}(t_n)\left(\frac{e^{t_n}}{e^{t_n}-1}\right)^{K-1} < e^{2e/e^{t_1}-1}\left(\frac{e^{t_1}}{e^{t_1}-1}\right)^{K-1},$$

the first inequality following from Lemma 4, the second following from the fact that both $C(t_n)$ and $(\frac{e^{t_n}}{e^{t_n}-1})^{K-1}$ achieve their maximum values at $n=1$.

As in Lemma 3, we have $\mathsf{E}[\psi_{mK}^4(S_m)] = \exp(6K^2 t_m)$ and $\mathsf{E}[\psi_{mK}^2(S_m)] = \exp(K^2 t_m)$. Making these substitution in $A_K$ and in (39), we get

$$\begin{aligned}(42)\quad \mathsf{E}[\|\hat{C}_n - C_n\|_n^2] &< (2^{m-n} - 1)\frac{(K+1)^{m-n+2}}{N} \\ &\quad \times B_K H_K^{m-n} e^{6(m-n)K^2 t_m + 2K^2 t_m}(1 + o(1)).\end{aligned}$$

The factor $(K+1)^2$ is negligible compared to the exponential factor in (42). The factors $B_K$ and $H_K$ grow exponentially in $K$, but their exponents are linear in $K$, whereas the dominant exponent in (42) is quadratic in $K$. Thus, $B_K$ and $H_K$ are also negligible for large $K$. If we set

$$K = \sqrt{\frac{(1-\delta)\log N}{(6(m-n)+2)t_m}}$$

for any $\delta > 0$, then

$$\lim_{N\to\infty} \sup_{|\beta_{m-1}|=1} \mathsf{E}[\|\hat{C}_n - C_n\|_n^2] = 0.$$

On the other hand, we know from Theorem 2 that if

$$K = \sqrt{\frac{(1+\delta)\log N}{3t_m + t_{m-1}}}$$

for any $\delta > 0$, then

$$\lim_{N\to\infty} \sup_{|\beta_{m-1}|=1} \mathsf{E}[\|\hat{C}_{m-1} - C_{m-1}\|_{m-1}^2] = \infty.$$

Thus, the crititcal rate of $K$ for the multiperiod problem is $O(\sqrt{\log N})$, just as in the single-period problem.

## 5. Proofs for the single-period problem.

5.1. *Normal setting.*

PROOF OF LEMMA 1. Equation (28) follows immediately from the orthogonality and martingale properties of the Hermite polynomials. Using

(24), we get

$$\mathsf{E}[\psi_{2k_2}^2(S_2)\psi_{1k_1}^2(S_1)] = \mathsf{E}\left[\left(\frac{He_{k_2}(S_2/\sqrt{t_2})}{\sqrt{k_2!}}\frac{He_{k_1}(S_1/\sqrt{t_1})}{\sqrt{k_1!}}\right)^2\right]$$

$$= (k_1!k_2!)\mathsf{E}\left[\sum_{k=0}^{k_2}\frac{He_{2k}(S_2/\sqrt{t_2})}{(k!)^2(k_2-k)!}\sum_{l=0}^{k_1}\frac{He_{2l}(S_1/\sqrt{t_1})}{(l!)^2(k_1-l)!}\right]$$

$$= (k_1!k_2!)\sum_{k=0}^{k_2}\sum_{l=0}^{k_1}\frac{\mathsf{E}[He_{2k}(S_2/\sqrt{t_2})He_{2l}(S_1/\sqrt{t_1})]}{(k!)^2(k_2-k)!(l!)^2(k_1-l)!}$$

$$= (k_1!k_2!)\sum_{k=0}^{k_1\wedge k_2}\frac{(2k)!\,(t_1/t_2)^k}{(k!)^2(k_2-k)!(k!)^2(k_1-k)!}$$

$$= \sum_{k=0}^{k_1\wedge k_2}\rho^{-k}\binom{2k}{k}\binom{k_1}{k}\binom{k_2}{k}.$$

The fourth equality applies (28). □

PROOF OF LEMMA 2. The ratio between the $(k+1)$st summand and the $k$th summand in (30) is

$$r_{kK} = \frac{\rho^{-(k+1)}\binom{2k+2}{k+1}\binom{K}{k+1}^2}{\rho^{-k}\binom{2k}{k}\binom{K}{k}^2} = \frac{2(2k+1)}{\rho(k+1)}\frac{(K-k)^2}{(k+1)^2}.$$

For $0 \leq k \leq K-1$, its derivative with respect to $k$ is

$$\frac{1}{\rho(k+1)^4}(8kK(k-K) + 4(k-K) + 10k^2 - 8kK - 2K^2) < 0.$$

Thus, $r_{kK}$ is strictly decreasing in $k$. At $k=0$, $r_{kK} = 4K^2/\rho$, which is greater than 1 for all sufficiently large $K$; and at $k = K-1$,

$$r_{K-1,K} = \frac{2(2K-1)}{\rho K^3},$$

which is less than 1 for all $K \geq 2$. Thus, for all sufficiently large $K$, $k^*$ is characterized by the condition

$$k^* = \min\{k : r_{kK} \leq 1\}.$$

The condition $r_{kK} \leq 1$ is equivalent to

(43) $$\frac{4(K-k)^2}{\rho(k+1)^2} \leq \frac{2k+2}{2k+1},$$

and $(2k+2)/(2k+1)$ is greater than 1 for all positive $k$. The ratio on the left-hand side of (43) is decreasing in $k$, $0 \leq k \leq K-1$, so if we define

$$(44) \qquad k_1 = \min\left\{k : \frac{4(K-k)^2}{\rho(k+1)^2} \leq 1\right\},$$

then $k^* \leq k_1$.

For any fixed $k$, the inequality in (43) will be violated for all sufficiently large $K$, so $k^*$ must increase without bound as $K \to \infty$. It follows that $(2k^*+2)/(2k^*+1) \to 1$. If for some $\varepsilon > 0$, we define

$$k_2 = \min\left\{k : \frac{4(K-k)^2}{\rho(k+1)^2} \leq 1+\varepsilon\right\},$$

then $k^* \geq k_2$ for all sufficiently large $K$. Thus, $k_2 \leq k^* \leq k_1$.

For $k_1$, we examine the equation

$$\frac{4(K-k)^2}{\rho(k+1)^2} = 1.$$

The only root of this equation less than $K$ is

$$\hat{k} = \frac{2K - \sqrt{\rho}}{2 + \sqrt{\rho}} = \left(\frac{2}{2+\sqrt{\rho}}\right) K(1+o(1)).$$

The solution $k_1$ to (44) is either $\lfloor \hat{k} \rfloor$ or $\lfloor \hat{k} \rfloor + 1$, so $k_1/\hat{k} \to 1$.

The same argument applied to the equation

$$\frac{4(K-k)^2}{\rho(k+1)^2} = 1+\varepsilon$$

shows that

$$k_2 = \left(\frac{2}{2+\sqrt{\rho(1+\varepsilon)}}\right) K(1+o(1)).$$

Noting that we may take $\varepsilon > 0$ arbitrarily small and $k_2 \leq k^* \leq k_1$ concludes the proof. $\square$

### 5.2. Lognormal setting.

PROOF OF LEMMA 3. Using the martingale property of $\psi_k(S(t))$ and the moment generating function of $W(t_1)$, we get

$$\begin{aligned}
\mathsf{E}[\psi_{k_1}(S(t_1))\psi_{k_2}(S(t_2))] &= \mathsf{E}[\mathsf{E}[\psi_{k_1}(S(t_1))\psi_{k_2}(S(t_2))|W(t_1)]] \\
&= \mathsf{E}[\psi_{k_1}(S(t_1))\psi_{k_2}(S(t_1))] \\
&= \mathsf{E}[e^{(k_1+k_2)W(t_1) - k_1^2 t_1/2 - k_2^2 t_2/2}] \\
&= e^{k_1 k_2 t_1}.
\end{aligned}$$

The second part of the lemma works similarly. □

PROOF OF LEMMA 4. The first assertion follows from the observation that the largest entry of $\Psi$ is $e^{K^2 t}$. For the second assertion, we note that $\Psi$ has the form of a Vandermonde matrix, allowing calculation of its determinant (using [12], page 322),

$$\det \Psi = \prod_{0 \leq q < r \leq K} (e^{rt} - e^{qt}). \tag{45}$$

By standard linear algebra, the inverse of $\Psi$ is given by

$$\Psi^{-1} = \frac{\Psi^*}{\det \Psi}, \tag{46}$$

where

$$\Psi^*_{qr} = (-1)^{q+r} \det \Psi(q|r),$$

and $\Psi(q|r)$ denotes the matrix obtained by deleting the $q$th row and $r$th column from $\Psi$. Two cases arise, depending on whether $q = r = 1$ or not.

*Case* 1. $q \neq 1$ or $r \neq 1$. Since $\Psi$ is symmetric, $\det \Psi(q|r) = \det \Psi(r|q)$, so it suffices to suppose $r \neq 1$. We can then compute the determinant of $\Psi(q|r)$ using [12], page 333. Through (46) this leads to

$$\Psi^{-1}_{qr} = \frac{(-1)^{q+r} \sum_{s_1 < \cdots < s_{K-(r-1)}, s_d \neq q-1} \exp\{\sum_{d=1}^{K-(r-1)} s_d t\}}{\prod_{j=0}^{q-2} (e^{(q-1)t} - e^{jt}) \prod_{j=q}^{K} (e^{jt} - e^{(q-1)t})}, \tag{47}$$

the sum ranging over $s_1, \ldots, s_d$ taking values in $\{0, \ldots, K\}$.

The lemma requires an upper bound on the numerator and a lower bound on the denominator. To bound the numerator, for $\hat{r} = 1, \ldots, K-1$, set

$$R(K, q, \hat{r}) = \sum_{s_1 < \cdots < s_{\hat{r}}, s_d \neq q-1} \exp\left\{\sum_{d=1}^{\hat{r}} s_d t\right\}.$$

We now claim that

$$R(K, q, \hat{r}) < R(K, 1, \hat{r}) < \frac{e^{\hat{r}(K+1)t}}{(e^t - 1)^{\hat{r}} e^{\hat{r}(\hat{r}-1)t/2}} \tag{48}$$

for $\hat{r} = 1, \ldots, K-1$. That $R(K, q, \hat{r}) < R(K, 1, \hat{r})$ is immediate from the definition of $R(K, q, \hat{r})$. The second inequality is proved by induction in $\hat{r}$. When $\hat{r} = 1$,

$$R(K, 1, \hat{r}) = (e^t + \cdots + e^{Kt}) < \frac{e^{(K+1)t}}{e^t - 1}.$$

Then

$$R(K,1,\hat{r}+1) = \sum_{s_1<\cdots<s_{\hat{r}+1}, s_d\neq 0} \exp\left\{\sum_{d=1}^{\hat{r}+1} s_d t\right\}$$

$$= \sum_{i_1=1}^{K-\hat{r}} \sum_{i_2=i_1+1}^{K-\hat{r}+1} \cdots \sum_{i_{\hat{r}+1}=\hat{r}+1}^{K} \exp\left\{\sum_{j=1}^{\hat{r}+1} i_j t\right\}$$

$$= \sum_{i_1=1}^{K-\hat{r}} e^{i_1 t} R(K,1,\hat{r})$$

$$< \sum_{i_1=1}^{K-\hat{r}} e^{i_1 t} \frac{e^{\hat{r}(K+1)t}}{(e^t-1)^{\hat{r}} e^{\hat{r}(\hat{r}-1)t/2}}$$

$$< \frac{e^{(K+1-\hat{r})t}}{e^t-1} \frac{e^{\hat{r}(K+1)t}}{(e^t-1)^{\hat{r}} e^{\sum_{j=1}^{\hat{r}-1} jt}}$$

(49)
$$= \frac{e^{(\hat{r}+1)(K+1)t}}{(e^t-1)^{\hat{r}+1} e^{\hat{r}(\hat{r}+1)t/2}}.$$

Thus, (48) holds.

The fact that

$$\frac{\partial}{\partial \hat{r}}\left(\frac{e^{\hat{r}(K+1)t}}{(e^t-1)^{\hat{r}} e^{\hat{r}(\hat{r}-1)t/2}}\right) > 0,$$

implies that (49) achieves its maximum when $\hat{r} = K-1$. Thus,

(50) $$R(K,q,\hat{r}) < \frac{e^{(K-1)(K+1)t}}{(e^t-1)^{K-1} e^{(K-1)(K-2)t/2}},$$

for $q=1,\ldots,K+1$, $r=2,\ldots,K+1$.

Next, we show that the denominator of (47) is bounded below by $C(t)\exp(K \times (K+1)t/2)$, with $C(t) = \exp(-2e/(e^t-1)^2)$. For this, we rewrite the denominator of (47) as

$$\prod_{j=0}^{q-2}(e^{(q-1)t} - e^{jt}) \prod_{j=q}^{K}(e^{jt} - e^{(q-1)t})$$

$$= \prod_{j=0}^{q-2} e^{(q-1)t}\left(1 - \frac{e^{jt}}{e^{(q-1)t}}\right) \prod_{j=q}^{K} e^{jt}\left(1 - \frac{e^{(q-1)t}}{e^{jt}}\right)$$

(51) $$= e^{(q-1)^2 t + \sum_{j=q}^{K} jt} \prod_{j=1}^{q-1}\left(1 - \frac{1}{e^{jt}}\right) \prod_{j=1}^{K-q+1}\left(1 - \frac{1}{e^{(q-1)t}}\right)$$

$$> e^{K(K+1)t/2 + (q-1)(q-2)/2} \prod_{j=1}^{K}\left(1 - \frac{1}{e^{jt}}\right) \prod_{j=1}^{K}\left(1 - \frac{1}{e^{jt}}\right)$$

$$> e^{K(K+1)t/2} \prod_{j=1}^{K}\left(1 - \frac{1}{e^{jt}}\right)^2.$$

Taking the logarithm of the product over $j$ and applying a Taylor expansion yields terms of the form

$$\log\left(1 - \frac{1}{e^{jt}}\right) = -\frac{1}{e^{jt}} - \frac{1}{2e^{2jt}} - \cdots - \frac{1}{ne^{njt}} - \cdots$$

$$> -\frac{1}{e^{jt}} - \frac{1}{e^{(j+1)t}} - \cdots - \frac{1}{e^{(n-1+j)t}} - \cdots$$

$$= -\frac{1}{e^{jt}} \frac{e^t}{e^t - 1}.$$

Therefore,

$$\sum_{j=1}^{K} \log\left(1 - \frac{1}{e^{jt}}\right) > -\frac{e^t}{e^t - 1} \sum_{j=1}^{K} \frac{1}{e^{jt}} = \frac{-e^t}{(e^t - 1)^2}(1 + o(1))$$

and

(52) $$\prod_{j=1}^{K}\left(1 - \frac{1}{e^{jt}}\right) > e^{-e/(e^t - 1)^2}.$$

Finally, by (51) and (52), we get that the denominator of (47) is bounded below by

(53) $$e^{-2e/(e^t - 1)^2} e^{K(K+1)t/2}.$$

Applying this lower bound and (50) to (47), we get

(54) $$\Psi_{qr}^{-1} < e^{2e/(e^t - 1)^2}\left(\frac{e^t}{e^t - 1}\right)^{K-1} = C^{-1}(t)\left(\frac{e^t}{e^t - 1}\right)^{K-1}$$

for $q, r$ not both equal to 1.

*Case* 2. $q = 1$ and $r = 1$. Because $\Psi\Psi^{-1} = I$ and all entries of the first row of $\Psi$ are 1, we have

(55) $$|\Psi_{11}^{-1}| = \left|1 - \sum_{r=2}^{K+1} \Psi_{1r}^{-1}\right| < 1 + \sum_{r=2}^{K+1} |\Psi_{1r}^{-1}| < C^{-1}(t)K\left(\frac{e^t}{e^t - 1}\right)^K.$$

Combining (54) and (55) we get

$$\|\Psi^{-1}\| = \sqrt{\sum_{q,r}(\Psi_{qr}^{-1})^2} < C^{-1}(t)(K+1)K\left(\frac{e^t}{e^t - 1}\right)^{K-1}.$$

□

**6. Proofs for the multiperiod problem.** As a tool for proving Theorem 3, we introduce a second sequence of coefficient estimates $\tilde{\beta}_n$ and $\tilde{\gamma}_n$. At each $n$, $\tilde{\beta}_n$ is the vector of coefficients that would be obtained using the algorithm of Section 2.2 if the coefficients $\beta_{n+1}$ were known exactly. More explicitly,

$$\tilde{\beta}_n = \Psi_n^{-1}\left(\frac{1}{N}\sum_{i=1}^N V_{n+1}(S_{n+1}^{(i)})\psi_n(S_n^{(i)})\right)$$
$$\equiv \Psi_n^{-1}\tilde{\gamma}_n,$$

with $V_{n+1}$ as in (10). The distinction between this and Step 2 of the algorithm is that here $V_{n+1}$ uses the true coefficients $\beta_n$ [as in (9)], whereas $\hat{V}_{n+1}$ in (13) uses the estimated coefficients $\hat{\beta}_{n+1}$. The estimates $\tilde{\beta}_n$ and $\tilde{\gamma}_n$ are not computable in practice and are simply used as a device for the proof. From the coefficients $\tilde{\beta}_n$ define

$$\tilde{C}_n(x) = \sum_{k=0}^K \tilde{\beta}_{nk}\psi_{nk}(x) = (\hat{L}_n C_{n+1})(x).$$

Thus, $\tilde{C}_n$ results from applying the estimated operator $\hat{L}_n$ to the exact function $C_{n+1}$, whereas $\hat{C}_n$ results from applying the estimated operator to the estimated function $\hat{C}_{n+1}$.

The proof of Theorem 3 also relies on two lemmas.

LEMMA 5. *Under conditions* (A0) *and* (B1)–(B3),

$$|\gamma_{m-n}|^2 \leq (2H_K \mathsf{E}[\psi_{mK}^4(S_m)])^{n-1}(K+1)^{n+1}(\mathsf{E}[\psi_{mK}^2(S_m)])^2(1+o(1))$$

*for* $n = 1, \ldots, m-1$.

PROOF. First note that for any $x \in \Re$,

$$C_n^2(x) = (\psi_n^\top(x)\Psi_n^{-1}\gamma_n)^2 \leq |\psi_n(x)|^2 \|\Psi_n^{-1}\|^2 |\gamma_n|^2.$$

By the definition of $\gamma$, together with the fact $|\max\{a,b\}| \leq |a| + |b|$, we get

$$|\gamma_{n,k}| = |\mathsf{E}[\psi_{nk}(S_n)\max\{h_{n+1}(S_{n+1}), C_{n+1}(S_{n+1})\}]|$$
$$\leq \mathsf{E}[|\psi_{nk}(S_n)h_{n+1}(S_{n+1})|] + \mathsf{E}[|\psi_{nk}(S_n)\psi_{n+1}^\top(S_{n+1})\Psi_{n+1}^{-1}\gamma_{n+1}|]$$
$$\leq \sqrt{\mathsf{E}[\psi_{nk}^2(S_n)]\mathsf{E}[h_{n+1}^2(S_{n+1})]}$$
$$\quad + \|\Psi_{n+1}^{-1}\||\gamma_{n+1}|\sqrt{\mathsf{E}[\psi_{nk}^2(S_n)|\psi_{n+1}(S_{n+1})|^2]}$$
$$\leq \sqrt{c^K \mathsf{E}[\psi_{mK}^4(S_m)]} + B_K|\gamma_{n+1}|\sqrt{(K+1)}\sqrt{\mathsf{E}[\psi_{mK}^4(S_m)]}.$$

The last inequality uses (B1), (B3) and the inequality $\mathsf{E}[h^2] \le \sqrt{\mathsf{E}[h^4]}$. Thus,

$$|\gamma_n|^2 = \sum_{k=0}^{K} \gamma_{n,k}^2$$

(56)
$$\le 2(K+1)\mathsf{E}[\psi_{mK}^4(S_m)](c^K + B_K^2(K+1)|\gamma_{n+1}|^2)$$
$$\le 2(K+1)\mathsf{E}[\psi_{mK}^4(S_m)]H_K(1 + |\gamma_{n+1}|^2),$$

with $H_K = \max\{c^K, B_K^2(K+1)\}$ as defined in Section 4.1.

Conditions (B1) and (B2) imply that

$$|\gamma_{m-1}|^2 \le \|\Psi_{m-1}\|^2 \le (K+1)^2(\mathsf{E}[\psi_{mK}^2(S_m)])^2.$$

Then (56) gives

$$|\gamma_{m-2}|^2 \le 2(K+1)\mathsf{E}[\psi_{mK}^4(S_m)]H_K(1+|\gamma_{m-1}|^2)$$
$$= 2H_K\mathsf{E}[\psi_{mK}^4(S_m)](K+1)^3(\mathsf{E}[\psi_{mK}^2(S_m)])^2(1+o(1)),$$
$$|\gamma_{m-3}|^2 \le 2(K+1)\mathsf{E}[\psi_{mK}^4(S_m)]H_K(1+|\gamma_{m-2}|^2)$$
$$= (2H_K\mathsf{E}[\psi_{mK}^4(S_m)])^2(K+1)^4(\mathsf{E}[\psi_{mK}^2(S_m)])^2(1+o(1))$$

and, proceeding by induction, completes the proof. □

LEMMA 6. *Under conditions* (A0) *and* (B1)–(B3),

$$\mathsf{E}[\|\hat{C}_n - C_n\|_n^2] \le B_K \sum_{l=1}^{m-n} A_K^{m-n-l} \mathsf{E}[|\tilde{\gamma}_{m-l} - \gamma_{m-l}|^2].$$

PROOF. By the definition of $C, \hat{C}$ and $\tilde{C}$ and the triangle inequality, we have

$$\mathsf{E}[\|\hat{C}_n - C_n\|_n^2] = \mathsf{E}[\|\hat{L}_n\hat{C}_{n+1} - L_nC_{n+1}\|_n^2]$$
$$\le \mathsf{E}[\|\hat{L}_n\hat{C}_{n+1} - \hat{L}_nC_{n+1}\|_n^2 + \|\hat{L}_nC_{n+1} - L_nC_{n+1}\|_n^2].$$

Now,

$$\hat{L}_n\hat{C}_{n+1} - \hat{L}_nC_{n+1} = \psi_n^\top \Psi_n^{-1}(\hat{\gamma}_n - \tilde{\gamma}_n),$$

so

$$\|\hat{L}_n\hat{C}_{n+1} - \hat{L}_nC_{n+1}\|_n^2$$
$$= (\hat{\gamma}_n - \tilde{\gamma}_n)^\top \Psi_n^{-1}\left(\int \psi_n(x)\psi_n(x)^\top g_n(x)\,dx\right)\Psi_n^{-1}(\hat{\gamma}_n - \tilde{\gamma}_n)$$
$$= (\hat{\gamma}_n - \tilde{\gamma}_n)^\top \Psi_n^{-1}(\hat{\gamma}_n - \tilde{\gamma}_n)$$
$$\le \|\Psi_n^{-1}\||(\hat{\gamma}_n - \tilde{\gamma}_n)|^2.$$

The same bound holds with $\hat{L}_n C_{n+1}$ replaced by $L_n C_{n+1}$ and $\hat{\gamma}_n$ replaced by $\gamma_n$. Thus,

(57) $\qquad \mathsf{E}[\|\hat{C}_n - C_n\|_n^2] \leq B_K(\mathsf{E}[|\hat{\gamma}_n - \tilde{\gamma}_n|^2] + \mathsf{E}[|\tilde{\gamma}_n - \gamma_n|^2]).$

Using the definitions of $\hat{\gamma}_n$ and $\tilde{\gamma}_n$ and the inequality $|\max\{a,b\} - \max\{a,c\}| \leq |b-c|$, we get

(58)
$$(\hat{\gamma}_{nk} - \tilde{\gamma}_{nk})^2 \leq \left(\frac{1}{N}\sum_{i=1}^{N}|\psi_{nk}(S_n^{(i)})|\,|\max\{h_{n+1}(S_{n+1}^{(i)}), \hat{C}_{n+1}(S_{n+1}^{(i)})\} - \max\{h_{n+1}(S_{n+1}^{(i)}), C_{n+1}(S_{n+1}^{(i)})\}|\right)^2$$
$$\leq \left(\frac{1}{N}\sum_{i=1}^{N}|\psi_{nk}(S_n^{(i)})|\,|\hat{C}_{n+1}(S_{n+1}^{(i)}) - C_{n+1}(S_{n+1}^{(i)})|\right)^2$$
$$\leq \frac{1}{N}\sum_{i=1}^{N}\psi_{nk}^2(S_n^{(i)})(\hat{C}_{n+1}(S_{n+1}^{(i)}) - C_{n+1}(S_{n+1}^{(i)}))^2.$$

The paths $S^{(i)}$, $i = 1, \ldots, N$, in this expression are independent of the coefficients of $\hat{C}_{n+1}$ (see Step 2 of the algorithm), so

(59) $\qquad \mathsf{E}[(\hat{\gamma}_{nk} - \tilde{\gamma}_{nk})^2] = \mathsf{E}[\psi_{nk}^2(S_n)(\hat{C}_{n+1}(S_{n+1}) - C_{n+1}(S_{n+1}))^2],$

with $(S_n, S_{n+1})$ independent of the coefficients of $\hat{C}_{n+1}$.

To bound (59), we use
$$(\hat{C}_{n+1}(S_{n+1}) - C_{n+1}(S_{n+1}))^2 = (\psi_{n+1}^\top(S_{n+1})\Psi_{n+1}^{-1}(\hat{\gamma}_{n+1} - \gamma_{n+1}))^2$$
$$\leq |\psi_{n+1}^\top(S_{n+1})|^2\|\Psi_{n+1}^{-1}\|^2|\hat{\gamma}_{n+1} - \gamma_{n+1}|^2.$$

The independence of $(S_n, S_{n+1})$ and $\hat{\gamma}_{n+1}$ then gives

$$\mathsf{E}[\psi_{nk}^2(S_n)(\hat{C}_{n+1}(S_{n+1}) - C_{n+1}(S_{n+1}))^2]$$
$$\leq \|\Psi_{n+1}^{-1}\|^2 \mathsf{E}[\psi_{nk}^2(S_n)|\psi_{n+1}(S_{n+1})|^2]\mathsf{E}[|\hat{\gamma}_{n+1} - \gamma_{n+1}|^2]$$
$$\leq B_K^2(K+1)\mathsf{E}[\psi_{nK}^2(S_n)\psi_{n+1,K}^2(S_{n+1})]\mathsf{E}[|\hat{\gamma}_{n+1} - \gamma_{n+1}|^2]$$
$$\leq B_K^2(K+1)\sqrt{\mathsf{E}[\psi_{nK}^4(S_n)]\mathsf{E}[\psi_{n+1,K}^4(S_{n+1})]}\mathsf{E}[|\hat{\gamma}_{n+1} - \gamma_{n+1}|^2]$$
$$\leq B_K^2(K+1)\mathsf{E}[\psi_{mK}^4(S_m)]\mathsf{E}[|\hat{\gamma}_{n+1} - \gamma_{n+1}|^2],$$

the last inequality following from (B1). Using this bound with (58) and (59), we get

$$\mathsf{E}[|\hat{\gamma}_n - \tilde{\gamma}_n|^2] = \sum_{k=0}^{K}\mathsf{E}[(\hat{\gamma}_{n,k} - \tilde{\gamma}_{n,k})^2]$$

$$\leq (K+1)^2 B_K^2 \mathsf{E}[\psi_{mK}^4(S_m)]\mathsf{E}[|\hat{\gamma}_{n+1} - \gamma_{n+1}|^2]$$

(60)
$$\leq A_K \mathsf{E}[|\hat{\gamma}_{n+1} - \gamma_{n+1}|^2]$$

(61)
$$\leq A_K \mathsf{E}[|\hat{\gamma}_{n+1} - \tilde{\gamma}_{n+1}|^2] + A_K \mathsf{E}[|\tilde{\gamma}_{n+1} - \gamma_{n+1}|^2].$$

By iteratively using (60)–(61), we get

$$\mathsf{E}[|\hat{\gamma}_n - \tilde{\gamma}_n|^2]$$
$$\leq A_K^{m-n-1}\mathsf{E}[|\hat{\gamma}_{m-1} - \gamma_{m-1}|^2] + \sum_{l=2}^{m-n-1} A_K^{m-n-l}\mathsf{E}[|\tilde{\gamma}_{m-l} - \gamma_{m-l}|^2]$$
$$= A_K^{m-n-1}\mathsf{E}[|\tilde{\gamma}_{m-1} - \gamma_{m-1}|^2] + \sum_{l=2}^{m-n-1} A_K^{m-n-l}\mathsf{E}[|\tilde{\gamma}_{m-l} - \gamma_{m-l}|^2]$$
$$= \sum_{l=1}^{m-n-1} A_K^{m-n-l}\mathsf{E}[|\tilde{\gamma}_{m-l} - \gamma_{m-l}|^2],$$

because $\hat{\gamma}_{m-1} = \tilde{\gamma}_{m-1}$ (since $\hat{C}_m = C_m = 0$). Using this bound in (57) concludes the proof. $\square$

PROOF OF THEOREM 3. Because each $\tilde{\gamma}_{nk}$ is an unbiased estimate of the corresponding $\gamma_{nk}$, $\mathsf{E}[(\tilde{\gamma}_{nk} - \gamma_{nk})^2]$ is the variance of $\tilde{\gamma}_{nk}$ and is therefore bounded above by the second moment of $\tilde{\gamma}_{nk}$. Thus,

$$\mathsf{E}[|\tilde{\gamma}_{m-n} - \gamma_{m-n}|^2]$$
$$= \sum_{k=0}^{K} \mathsf{E}[(\tilde{\gamma}_{m-n,k} - \gamma_{m-n,k})^2]$$
$$\leq \sum_{k=0}^{K} \frac{1}{N}\mathsf{E}[\psi_{m-n,k}^2(S_{m-n})\max\{h_{m-n+1}^2(S_{m-n+1}), C_{m-n+1}^2(S_{m-n+1})\}]$$
$$\leq \sum_{k=0}^{K} \frac{1}{N}\mathsf{E}[\psi_{m-n,k}^2(S_{m-n})(h_{m-n+1}^2(S_{m-n+1}) + C_{m-n+1}^2(S_{m-n+1}))]$$
$$\leq \sum_{k=0}^{K} \frac{1}{N}\mathsf{E}[\psi_{m-n,k}^2(S_{m-n})h_{m-n+1}^2(S_{m-n+1})]$$

(62)
$$+ \sum_{k=0}^{K} \frac{1}{N}\mathsf{E}[\psi_{m-n,k}^2(S_{m-n})\|\Psi_{m-n}^{-1}\|^2|\gamma_{m-n}|^2|\psi_{m-n+1}(S_{m-n+1})|^2].$$

AMERICAN OPTION PRICING 29For the first term in (62) we use the Cauchy–Schwarz inequality, (B1) and (B3) to get

$$\sum_{k=0}^{K} \frac{1}{N} \mathsf{E}[\psi_{m-n,k}^2(S_{m-n}) h_{m-n+1}^2(S_{m-n+1})]$$

(63)
$$\leq \frac{K+1}{N} \sqrt{\mathsf{E}[\psi_{mK}^4(S_m)] \mathsf{E}[h_{m-n+1}^4(S_{m-n+1})]}$$

$$\leq \frac{K+1}{N} H_K \mathsf{E}[\psi_{mK}^4(S_m)].$$

For the second term in (62) we again use Cauchy–Schwarz and (B1) to get

$$\sum_{k=0}^{K} \frac{1}{N} \mathsf{E}[\psi_{m-n,k}^2(S_{m-n}) \|\Psi_{m-n}^{-1}\|^2 |\gamma_{m-n}|^2 |\psi_{m-n+1}(S_{m-n+1})|^2]$$

(64)
$$\leq \frac{K+1}{N} H_K \mathsf{E}[\psi_{mK}^4(S_m)] |\gamma_{m-n}|^2.$$

Combining (62)–(64) and Lemma 5 we arrive at

$$\mathsf{E}[|\tilde{\gamma}_{m-n} - \gamma_{m-n}|^2]$$

$$\leq \frac{(K+1)^{n+2}}{N} 2^{n-1} (H_K \mathsf{E}[\psi_{mK}^4(S_m)])^n (\mathsf{E}[\psi_{mK}^2(S_m)])^2 (1+o(1))$$

$$= \frac{2^{n-1}(K+1)^2}{N} A_K^n (\mathsf{E}[\psi_{mK}^2(S_m)])^2 (1+o(1)).$$

By Lemma 6, we now get

$$\mathsf{E}[\|\hat{C}_n - C_n\|_n^2]$$

$$\leq B_K \left( \sum_{l=1}^{m-n} A_K^{m-n-l} \mathsf{E}[|\tilde{\gamma}_{m-l} - \gamma_{m-l}|^2] \right)$$

$$\leq B_K \frac{(K+1)^2}{N} A_K^{m-n} (\mathsf{E}[\psi_{mK}^2(S_m)])^2$$

$$\times (1 + 2 + \cdots + 2^{m-n-1})(1+o(1))$$

$$= (2^{m-n} - 1) B_K \frac{(K+1)^2}{N} A_K^{m-n} (\mathsf{E}[\psi_{mK}^2(S_m)])^2 (1+o(1)),$$

which concludes the proof. □

**7. Concluding remarks.** It is natural to ask to what extent our results depend on the fact that the basis functions we consider are polynomials. Some insight into this question can be gleaned from the analysis of the lower bound on $\mathsf{MSE}(\tilde{\beta})$ in the proof of Theorem 1. The lower bound results from choosing $Y = a_K \psi_{2K}(S_2)$ and its growth is driven by the second moment $a_K^2 \mathsf{E}[\psi_{2K}^2(S_2)\psi_{1K}^2(S_1)]$. With $\psi_{1K}$ orthogonal to the other basis functions at $t_1$, the condition $|\beta| = 1$ translates to $a_K = 1/\mathsf{E}[\psi_{2K}(S_2)\psi_{1K}(S_1)]$. Thus, the growth of the lower bound is driven by the growth of the ratio

$$\frac{\mathsf{E}[\psi_{2K}^2(S_2)\psi_{1K}^2(S_1)]}{(\mathsf{E}[\psi_{2K}(S_2)\psi_{1K}(S_1)])^2}$$

as $K$ increases. A few examples show that this ratio does indeed grow with $K$ even for choices of functions that grow much less quickly than polynomials. In the case of Brownian motion, explicit calculations show that for $\psi_{jK}(x) = \mathbf{1}\{x > K\}$, the ratio is $O(K \exp(K^2/2t_1))$ and for $\psi_{jK}(x) = \max\{0, x - K\}$, the ratio is $O(K^3 \exp(K^2/2t_1))$, so in both of these cases the growth rate is even faster than for the polynomials in Theorem 1. With $\psi_{jK}(x) = x^K \exp(-x)$, numerical calculations indicate that the ratio is roughly linear in $K$ (thus requiring roughly linear growth of $N$), but its magnitude is very large even at small values of $K$. These simple illustrations suggest that the phenomena observed in this paper may occur more generally. But see [7] for more positive results using bounded basis functions.

**Acknowledgment.** We thank the referee for a careful reading of the manuscript and helpful comments and corrections.

GRADUATE SCHOOL OF BUSINESS
COLUMBIA UNIVERSITY
NEW YORK, NEW YORK 10027
USA
E-MAIL: pg20@columbia.edu
E-MAIL: by52@columbia.edu